# SEQUENTIAL IMPORTANCE SAMPLING FOR MULTIWAY TABLES[1]

By Yuguo Chen, Ian H. Dinwoodie and Seth Sullivant

*University of Illinois at Urbana-Champaign, Duke University and Harvard University*

We describe an algorithm for the sequential sampling of entries in multiway contingency tables with given constraints. The algorithm can be used for computations in exact conditional inference. To justify the algorithm, a theory relates sampling values at each step to properties of the associated toric ideal using computational commutative algebra. In particular, the property of interval cell counts at each step is related to exponents on lead indeterminates of a lexicographic Gröbner basis. Also, the approximation of integer programming by linear programming for sampling is related to initial terms of a toric ideal. We apply the algorithm to examples of contingency tables which appear in the social and medical sciences. The numerical results demonstrate that the theory is applicable and that the algorithm performs well.

**1. Introduction.** Sampling from multiway contingency tables with given constraints can be used to compute exact Monte Carlo $p$-values of goodness-of-fit and parameter significance for conditional inference. This is desirable when the tables of interest are numerous but have entries that raise doubts about the validity of asymptotic methods. A classical application is testing for Hardy–Weinberg equilibrium with multiple alleles, where some alleles may be quite rare and result in sparse tables [20]. Other applications are described in [2, 6, 13]. A more general problem is sampling from nonnegative integer lattice points. This includes contingency tables, and further applications such as Monte Carlo EM algorithms with incomplete data [31] and Bayesian computation of posterior distributions [30].

Received July 2004; revised February 2005.
[1]Supported by NSF Grants DMS-02-00888, DMS-02-03762 and DMS-05-03981 and NSF Grant DMS-01-12069 to SAMSI.

Markov chain Monte Carlo (MCMC) has been a popular technique for generating random samples from tables with given constraints. It is usually easy to program, does not require a lot of memory, and has wide applicability. Diaconis and Sturmfels [14] gave algebraic characterizations of the moves necessary to run such a Markov chain. However, for some loglinear models the constraints from sufficient statistics on multiway tables make it difficult to design irreducible Markov chains. Diaconis and Sturmfels [14] gave a method to produce Markov moves that connect all tables with given constraints, but in some practical cases, such as large logistic regression examples, the moves cannot be computed. It is sometimes possible to do computations with a smaller collection of moves by letting some entries in the space of tables go negative. This idea is used in [4, 7]. The cost is a longer running time for the Markov chain. In general, the running times of these Markov chains are very difficult to judge. Therefore, Markov chains have three disadvantages: (1) they can be hard to design, (2) they can take a long time to run to stationarity, and (3) the time to run to stationarity may not be clear.

Sequential importance sampling (SIS) avoids these disadvantages of Markov chains because it is relatively easy to implement and there is no issue of converging to a stationary distribution. Chen et al. [6] introduced an SIS procedure for simulating two-way zero–one and contingency tables with fixed marginal sums, which compares favorably with other existing Monte Carlo-based algorithms. Similar techniques have also been applied to a logistic regression problem in [7]. This paper shows that SIS can be implemented efficiently for many multiway contingency table problems that have been studied mostly with Markov chains.

The idea behind SIS is to sample cell entries in the contingency table one after the other so that the final joint distribution (i.e., the *proposal distribution*) is close to the target distribution. SIS does not have the same disadvantages as a Markov chain, because the method terminates at the last cell and generates i.i.d. samples from the proposal distribution. However, SIS raises a new set of implementation issues. The main problems are approximating the support of the marginal distribution of each cell quickly, and then approximating the marginal distribution on the support set with a proposal distribution. We show how properties of the sampling set at each step can be deduced from algebraic conditions on a collection of Markov moves. The results of this paper extend the applicability of SIS from two-way tables [6] to a wider range of multiway tables and allow further comparison with Markov chain methods.

The target distribution on the collection of tables may be hypergeometric, which arises in conditional inference with multinomial sampling, or it may be another related distribution such as the one for Hardy–Weinberg proportions. SIS can yield an approximate count of constrained tables very quickly when the target distribution is uniform. This application has been carried out in [6], where SIS was shown to be more efficient than Markov chains for counting and testing two-way tables. Combinatorists are interested in counting tables with given constraints [11]. Counting tables is also related to conditional volume tests [13]. In our multiway examples, we found approximate counts of tables without difficulty. The exact counting software LattE [11] confirmed the counts on the two smaller examples. The uniform target distribution is also useful for Bayesian applications where a uniform prior on probabilities leads to equally likely tables, and for the conditional volume test [13].

The paper is organized as follows. In Section 2 we introduce essential ideas of SIS. The algebraic conditions for efficient sampling are formulated in Sections 3 and 4. Many of the algebraic ideas of Markov chains on lattice points are used. Section 3 treats the basic case where properties of polynomials generating the toric ideal are related to SIS. Section 4 is more technical and develops stronger methods for subsets of the Markov basis. These results can apply when the observed margins imply conditions of positivity on the tables constrained by the margin values.

Section 5 is about the relationship between linear programming (LP) and integer programming (IP). When the support of the marginal cell distribution is an interval of integers $[l, u]$, a situation established under conditions in Sections 3 and 4 and which occurs often in practice, one needs the values of the upper and lower bounds. Knowing then that LP and IP give nearly the same answer is important, because using an IP algorithm at each step in the procedure would be much slower than using LP. A precise algebraic relationship between LP and IP is developed in [22], which gives an algorithm for finding the maximum difference between the two over all conceivable data sets. The results here may be easier to apply in some examples. In practice it is not essential that LP and IP be identical. Section 6 discusses sampling distributions for different target distributions. In Section 7 we give a range of examples to show how well SIS can work in real problems. Section 8 provides concluding remarks.

**2. Elements of SIS.** Let $\Omega$ denote the set of all contingency tables with given constraints. Assume $\Omega$ is nonempty. The $p$-value for conditional inference on contingency tables can often be written as

$$\mu = E_p f(\mathbf{n}) = \sum_{\mathbf{n} \in \Omega} f(\mathbf{n}) p(\mathbf{n}), \tag{1}$$

where $p(\mathbf{n})$ is the underlying distribution on $\Omega$, which is usually uniform or hypergeometric and only known up to a normalizing constant, and $f(\mathbf{n})$ is a function of the test statistic. For example, if we let

$$f(\mathbf{n}) = \mathbb{1}_{\{p(\mathbf{n}) \leq p(\mathbf{n}_0)\}}, \tag{2}$$

where $\mathbf{n}_0$ is the observed table, formula (1) gives the $p$-value of the exact test [20]. In many cases sampling from $p(\mathbf{n})$ directly is difficult. The importance sampling approach is to simulate a table $\mathbf{n} \in \Omega$ from a different distribution $q(\cdot)$, where $q(\mathbf{n}) > 0$ for all $\mathbf{n} \in \Omega$, and estimate $\mu$ by

$$(3) \qquad \hat{\mu} = \frac{\sum_{i=1}^{N} f(\mathbf{n}_i) p(\mathbf{n}_i)/q(\mathbf{n}_i)}{\sum_{i=1}^{N} p(\mathbf{n}_i)/q(\mathbf{n}_i)},$$

where $\mathbf{n}_1, \ldots, \mathbf{n}_N$ are i.i.d. samples from $q(\mathbf{n})$. We can also estimate the total number of tables in $\Omega$ by

$$(4) \qquad \widehat{|\Omega|} = \frac{1}{N} \sum_{i=1}^{N} \frac{\mathbb{1}_{\{\mathbf{n}_i \in \Omega\}}}{q(\mathbf{n}_i)},$$

because $|\Omega| = \sum_{\mathbf{n} \in \Omega} \frac{1}{q(\mathbf{n})} q(\mathbf{n})$. The underlying distribution on $\Omega$ corresponding to this case is uniform.

In order to evaluate the efficiency of an importance sampling algorithm, we can look at the number of i.i.d. samples from the target distribution that are needed to give the same standard error for $\hat{\mu}$ as $N$ importance samples. A rough approximation for this number is the *effective sample size* [24]

$$(5) \qquad \text{ESS} = \frac{N}{1 + cv^2},$$

where the *coefficient of variation* ($cv$) is defined as

$$(6) \qquad cv^2 = \frac{\text{var}_q\{p(\mathbf{n})/q(\mathbf{n})\}}{E_q^2\{p(\mathbf{n})/q(\mathbf{n})\}}.$$

Accurate estimation generally requires a low $cv^2$, that is, $q(\mathbf{n})$ must be sufficiently close to $p(\mathbf{n})$. We will use $cv^2$ as a measure of efficiency for an importance sampling scheme. In practice, the theoretical value of $cv^2$ is unknown, so its sample counterpart is used to estimate $cv^2$. The standard error of $\hat{\mu}$ or $\widehat{|\Omega|}$ can be simply estimated by further repeated sampling [6].

SIS as it applies to multiway tables fills in the entries of a table cell by cell, in a way that guarantees that every table in $\Omega$ can be produced. More precisely, we stack all entries of the table into a long vector $\mathbf{n}$, and start by sampling the first cell count $n_1$ of the vector $\mathbf{n}$ with a proposal distribution $q(n_1)$. Conditional on the realization of the first cell, we sample the second cell count $n_2$ with a proposal distribution $q(n_2|n_1)$, and then move forward sequentially until all the cells are sampled. Denoting the cell counts of $\mathbf{n}$ by $n_1, \ldots, n_d$, we can write the joint proposal distribution $q$ as

$$q((n_1, \ldots, n_d)) = q(n_1) q(n_2|n_1) q(n_3|n_2, n_1) \cdots q(n_d|n_{d-1}, \ldots, n_1).$$

Ideally, one would like to sample a cell value from the marginal distribution of a cell entry, conditional on the entries that have already been sampled.

However, these marginal distributions are quite difficult to compute explicitly except in very small examples. SIS then raises some problems if it is to be used effectively: (1) When and how can the support of the marginal distribution $n_i|(n_{i-1},\ldots,n_1)$ be quickly determined or approximated? (2) How can the support of the marginal distribution be sampled with a proposal distribution $q$ that is close to the true underlying distribution $p$? We address these questions in the following sections of the paper.

**3. Sequential intervals and algebra.** When they apply SIS to the problem of sampling two-way contingency tables with fixed marginal sums, Chen et al. [6] notice that the support of the marginal distribution $n_i|(n_{i-1},\ldots,n_1)$ is an interval of integers [here $\mathbf{n} = (n_1,\ldots,n_d)$ is the table in a vector format]. Therefore, they can sample a value from the interval at each step and always produce a table in $\Omega$, that is, every table satisfies the constraints. This saves a lot of computing time compared to rejection sampling. Another advantage of having this interval property is that one can find a good proposal distribution $q(n_i|n_{i-1},\ldots,n_1)$ more easily than in the situation where there are gaps in the support set.

SIS tends to perform better when the sequential interval property holds, but for general constraints on multiway tables, it is not always true that one can fill in entries in sequence and expect the range of feasible values to be an interval of integers. Examples where the sequential interval property does not hold are very sparse logistic regression [7], many 3-way tables with certain margin constraints (see [12] for the full range of difficulties with 3-way tables) and some triangular tables of genotype data when cells are sampled in certain orders. Typically, there may be a problem if the moves of a Markov basis involve changes in some entry that are of size $\pm 2$ or larger. A precise condition is more complicated and weaker than "no moves of size greater than 1," and may depend on the margin values and the order of the sequential sampling. In this section we give the basic theorems that are not related to the actual values of the margin constraints. In the next section we strengthen the results.

Now we introduce notation for lattice points and the algebra of polynomials that will be used in our study of SIS. Let $A$ be an $r \times d$ matrix of nonnegative integers, denoted $Z_+$. In applications $d$ is the number of cells in the table, and $r$ is the number of parameters (not necessarily free) in an exponential family model. $A$ is often referred to as the constraint matrix and $r$ is the total number of constraints. We assume that a sum of some nonempty subset of the rows of $A$ is a strictly positive vector. In applications with multinomial sampling, this will be immediate because the sample size is fixed, so the constant vector of ones is a row or is in the row space of $A$. For $\mathbf{t} \in Z_+^r$, let

$$A^{-1}[\mathbf{t}] := \{\mathbf{n} \in Z_+^d : A\mathbf{n} = \mathbf{t}\}.$$

This is a collection of tables with linear constraints, that is, the set of nonnegative integer points inside a polytope. The linear constraint value $\mathbf{t}$ will sometimes informally be called a margin constraint. The value of $\mathbf{t}$ will typically be the sufficient statistics for a loglinear model. Our primary goal is to sample from $A^{-1}[\mathbf{t}]$.

Let us first recall the notion of a Markov move on $A^{-1}[\mathbf{t}]$. If $\mathbf{m} \in \texttt{ker}_Z(A)$ (the null space of $A$ in the integers), then $\mathbf{m}$ is a Markov move. With a collection of such moves, one can define a symmetric Markov chain on $A^{-1}[\mathbf{t}]$ by starting at an initial state $\mathbf{n} \in A^{-1}[\mathbf{t}]$, and then uniformly choosing one of the moves $\mathbf{m}$ and a sign on the move, and then moving to the new state $\mathbf{n} \pm \mathbf{m}$ if this new vector is nonnegative (i.e., every entry is nonnegative). A Markov basis $M_A$ for $A$ is a subset of $\texttt{ker}_Z(A)$ such that, for each pair of vectors $\mathbf{u}, \mathbf{v} \in Z_+^d$ with $A\mathbf{u} = A\mathbf{v}$, there is a sequence of vectors $\mathbf{m}_i \in M_A$, $i = 1, \ldots, l$, such that

$$\mathbf{u} = \mathbf{v} + \sum_{i=1}^{l} \mathbf{m}_i,$$

$$\mathbf{0} \leq \mathbf{v} + \sum_{i=1}^{j} \mathbf{m}_i, \qquad j = 1, \ldots, l.$$

That is, two nonnegative vectors with the same linear constraints can be connected with a sequence of increments from $M_A$ while always maintaining the linear constraints and the nonnegativity.

Define the polynomial ring $Q[x_1, \ldots, x_d]$ in indeterminates (polynomial variables) $x_1, \ldots, x_d$, one for each cell. Define the toric ideal

$$I_A := \langle \mathbf{x}^\mathbf{n} - \mathbf{x}^\mathbf{m} : A\mathbf{n} = A\mathbf{m} \rangle,$$

where $\mathbf{x}^\mathbf{n} := x_1^{n_1} x_2^{n_2} \cdots x_d^{n_d}$ is the usual monomial notation for a nonnegative integer vector of exponents $\mathbf{n} = (n_1, \ldots, n_d)$. The way to go between Markov moves and polynomials is simple: order and number the cells in the table, create an indeterminate (polynomial variable) for each cell in the table, and put the positive Markov move cell values on one monomial, put the negative values on another monomial, then form the difference. For example, the Markov move $(1, -1, -1, 1)'$ can be denoted as $x_1 x_4 - x_2 x_3$. The choice of cell ordering can be important, as in Example 7.5.

There are two fundamental algebraic ideas related to Markov bases. For $\mathbf{m} \in Z^d$, define $\mathbf{m}^+ = \max\{\mathbf{0}, \mathbf{m}\}$, $\mathbf{m}^- = \max\{\mathbf{0}, -\mathbf{m}\}$, so $\mathbf{m} = \mathbf{m}^+ - \mathbf{m}^-$. The first fundamental result, shown by Diaconis and Sturmfels ([14], Theorem 3.1), is that a finite generating set of binomials $\{\mathbf{x}^{\mathbf{m}_i^+} - \mathbf{x}^{\mathbf{m}_i^-}, i = 1, \ldots, g\}$ for $I_A$ defines Markov moves $\pm(\mathbf{m}_i^+ - \mathbf{m}_i^-)$, $i = 1, \ldots, g$, that are a Markov basis in that they connect all of $A^{-1}[\mathbf{t}]$ when chosen randomly as vector increments, regardless of the actual value of $\mathbf{t}$. In other words, a Markov basis always exists independently of the actual values of the linear constraints. The second fundamental result ([29], Theorem 8.14) is that a collection of moves will connect two tables **n** and **m** if $\mathbf{x^n} - \mathbf{x^m} \in I$, where $I$ is the ideal generated by the collection of moves. This is used to show connectivity for subcollections of the full Markov basis for particular values of **t** in Section 4.

DEFINITION 3.1. Define the projection operator $\pi_1 : Z^d \to Z$ by $\pi_1(z_1, \ldots, z_d) = z_1$.

LEMMA 3.1. *Suppose a Markov basis $M_A$ satisfies $\pi_1(M_A) \subset \{-1, 0, +1\}$. Then $\pi_1(A^{-1}[\mathbf{t}])$ is an interval of integers $[l_1, u_1]$.*

PROOF. One can connect tables $\mathbf{m}, \mathbf{n} \in A^{-1}[\mathbf{t}]$ with values $m_1$ and $n_1$ in the first coordinate by changing the first coordinate only $\pm 1$ at each step, so the gap between possible values cannot be greater than 1. □

If the columns of $A$ are $\mathbf{a}_1, \ldots, \mathbf{a}_d$, let $A_i = (\mathbf{a}_i, \mathbf{a}_{i+1}, \ldots, \mathbf{a}_d)$ be the matrix that deletes the first $i - 1$ columns and keeps the last $d - i + 1$ columns of $A$.

DEFINITION 3.2. The polytope $A^{-1}[\mathbf{t}]$ has the *sequential interval property* if $\pi_1(A^{-1}[\mathbf{t}])$ is an interval of integers $[l_1, u_1]$, and for $i = 1, \ldots, d-1$: if $n_i \in \pi_1(A_i^{-1}[\mathbf{t} - n_1\mathbf{a}_1 - \cdots - n_{i-1}\mathbf{a}_{i-1}])$, then $\pi_1(A_{i+1}^{-1}[\mathbf{t} - n_1\mathbf{a}_1 - \cdots - n_{i-1}\mathbf{a}_{i-1} - n_i\mathbf{a}_i])$ is also an interval of integers $[l_{i+1}, u_{i+1}]$.

The next result is the most basic connection between the sequential interval property and the exponents of a lex basis for the toric ideal. An important point is that the condition does not require that all exponents in the Markov basis have magnitude at most 1. Rather, it requires that the exponent be at most 1 on the indeterminate $x_i$ (square-free in $x_i$) on the moves that involve only the present and future cells $i, i+1, \ldots, d$ in the lex basis. This point is important for many examples, including $3 \times 3 \times 3$ tables with no-3-way interaction (Example 7.4).

With a particular cell order, the indeterminates are typically ordered $x_1 > x_2 > \cdots > x_d$, and then one can introduce term orders. We primarily use the lexicographic term order (lex order), which totally orders monomials (or, equivalently, their vector exponents corresponding to tables) by declaring $\mathbf{x^n} > \mathbf{x^m}$ if and only if the first entry from the left in $\mathbf{n} - \mathbf{m}$ is positive (or **n** is after **m** in the dictionary sense). Cox, Little and O'Shea ([10], page 52) explain term orders, including the grevlex order that we use in Section 5 where the indeterminates are taken in reverse order $x_d > x_{d-1} > \cdots > x_1$.

In the following, we use the term "Gröbner basis," which is a special generating set for an ideal ([10], page 74). Lex Gröbner basis (or lex basis) will mean Gröbner basis with respect to lexicographic term order ([10], page 54) and reduced Gröbner basis is a unique representation ([10], page 90).

PROPOSITION 3.1. *Suppose a Markov basis $M_A = \{\pm \mathbf{m}_1, \ldots, \pm \mathbf{m}_g\}$ has the property that $G := \{\mathbf{x}^{\mathbf{m}_i^+} - \mathbf{x}^{\mathbf{m}_i^-}, i = 1, \ldots, g\}$ is a lex Gröbner basis with ordering $x_1 > x_2 > \cdots > x_d$ on indeterminates and suppose the elements of $G \cap Q[x_i, \ldots, x_d]$ are square-free in $x_i$ for each $i$. Then $A^{-1}[\mathbf{t}]$ has the sequential interval property for all $\mathbf{t}$.*

PROOF. By the elimination theorem ([10], page 113), the lex basis $G$ has the property that $G \cap Q[x_i, \ldots, x_d]$ is a Gröbner basis for the ideal $I_{A_i} = \langle \mathbf{x}^{\mathbf{m}} - \mathbf{x}^{\mathbf{n}}, A_i \mathbf{m} = A_i \mathbf{n} \rangle$. Hence, by Theorem 3.1 of Diaconis and Sturmfels [14], the difference of the exponents (together with signs $\pm$) of elements in $G \cap Q[x_i, \ldots, x_d]$ is a Markov basis with 0 in coordinates $1, 2, \ldots, i-1$. An application of Lemma 3.1 to the matrix $A_i$ with first coordinate $n_i$ completes the proof. □

When using this result, some orders on the cells may have the square-free property and others may not, so it can be used to find good orderings on the cells. The sensitivity to cell ordering shows up in many examples, including logistic regression and Hardy–Weinberg testing with genotype data (Example 7.5).

In fact, the converse to Proposition 3.1 is also true, in the sense that matrices $A$, such that $A^{-1}[\mathbf{t}]$ has the sequential interval property regardless of $\mathbf{t}$, are characterized by their lex Gröbner bases.

PROPOSITION 3.2. *Let $A$ be a nonnegative integer matrix such that $A^{-1}[\mathbf{t}]$ has the sequential interval property for all $\mathbf{t}$. Then the reduced lex Gröbner basis $G$ for $I_A$ with ordering $x_1 > x_2 > \cdots > x_d$ has $G \cap Q[x_i, \ldots x_d]$ square-free in $x_i$ for all $i$.*

PROOF. It suffices to prove the claim on the first cell, the rest following by induction. Let $G := \{\mathbf{x}^{\mathbf{m}_i^+} - \mathbf{x}^{\mathbf{m}_i^-}\}$ be the reduced lex Gröbner basis. In particular, none of the monomials $\mathbf{x}^{\mathbf{m}_i^+}$ is divisible by the leading monomial of any other binomial in $I_A$. Suppose there is some $\mathbf{x}^{\mathbf{m}^+} - \mathbf{x}^{\mathbf{m}^-} \in G$ with $\pi_1(\mathbf{m}^+) = a > 1$. Let $\mathbf{t} = A\mathbf{m}^+$. Since $A^{-1}[\mathbf{t}]$ has the sequential interval property and $\pi_1(\mathbf{m}^-) = 0$, there exists $\mathbf{n} \in A^{-1}[\mathbf{t}]$ with $\pi_1(\mathbf{n}) = a - 1$. Then the binomial $x_1^{-a+1}(\mathbf{x}^{\mathbf{m}^+} - \mathbf{x}^{\mathbf{n}}) \in I_A$ is not equal to $\mathbf{x}^{\mathbf{m}^+} - \mathbf{x}^{\mathbf{m}^-}$, and has leading term $x_1^{-a+1}\mathbf{x}^{\mathbf{m}^+}$ which divides $\mathbf{x}^{\mathbf{m}^+}$. This is a contradiction and $\mathbf{x}^{\mathbf{m}^+} - \mathbf{x}^{\mathbf{m}^-}$ is not in the reduced Gröbner basis $G$. □

**4. Markov subbases.** In this section we give results that can be used when the full Markov basis does not have the required properties to guarantee sequential intervals. Situations where this occurs include logistic regression [7] and Example 7.3, where the lex bases for the toric ideals do not have the conditions of Proposition 3.1.

The results in this section use the particular values of the margin constraints, which may allow a smaller and simpler connecting set that we call a Markov subbasis. An existing method to study connectivity properties of subsets of a Markov basis is the primary decomposition ([10], page 208). While useful in some examples, it is usually quite difficult to compute. The methods in this section use computational tools that are more easily applied in many cases.

To motivate some of the ideas that follow, recall that in some contingency tables it is possible to easily identify a reasonable collection of Markov moves that preserve the required constraints and are a basis in the linear algebra sense for the kernel of the constraint matrix. However, a basis in the linear algebra sense does not always give a Markov basis—the Markov basis allows you to connect all tables while remaining nonnegative, a condition not guaranteed by the linear algebra basis. The smaller collection, while not a Markov basis, may connect tables with certain margin values while remaining nonnegative. The linear algebra basis can be enlarged to a Markov basis by a process called saturation discussed below, and the result can be much more complicated than the original collection of moves.

A lex basis for the toric ideal $I_A$ for a constraint matrix $A$ is quite special in that the Markov moves that involve cells $i, i+1, \ldots, d$ are a lex basis for the toric ideal for $I_{A_i}$. This is a consequence of the elimination theorem, and means that one lex basis calculation gives sequential sampling information about all the cells in sequence. With a collection of moves smaller than a lex basis for $I_A$, the theory is more difficult.

A Markov subbasis $M_{A,\mathbf{t}}$ for $\mathbf{t} \in Z_+^r$ and integer matrix $A$ is a finite subset of $\texttt{ker}_Z(A)$ such that, for each pair of vectors $\mathbf{u}, \mathbf{v} \in A^{-1}[\mathbf{t}]$, there is a sequence of vectors $\mathbf{m}_i \in M_{A,\mathbf{t}}$, $i = 1, \ldots, l$, such that

$$\mathbf{u} = \mathbf{v} + \sum_{i=1}^{l} \mathbf{m}_i,$$

$$\mathbf{0} \leq \mathbf{v} + \sum_{i=1}^{j} \mathbf{m}_i, \qquad j = 1, \ldots, l.$$

The connectivity through nonnegative lattice points only is required to hold for this specific $\mathbf{t}$.

LEMMA 4.1. *Suppose a Markov subbasis $M_{A,\mathbf{t}}$ satisfies $\pi_1(M_{A,\mathbf{t}}) \subset \{-1, 0, +1\}$. Then $\pi_1(A^{-1}[\mathbf{t}])$ is an interval of integers $[l_1, u_1]$.*

PROOF. One can connect tables with feasible values $n_1$ and $m_1$ in the first coordinate by changing the first coordinate only $\pm 1$ at each step, so the gap between possible values cannot be greater than 1. □

The following proposition is used in Examples 7.3 and 7.4, where Proposition 3.1 cannot be used. Recall that a lex basis for a toric ideal has the property that each elimination ideal ([10], page 113) is also a lex basis for a remaining toric ideal, so applying Lemma 4.1 in sequence is immediate. With a subbasis, however, one must add a technical condition involving saturation to get the sequential application of Lemma 4.1.

Saturation (see [28], page 113 or [25], page 215) is an algebraic procedure that enlarges an ideal. In our case the ideal will correspond to a collection of Markov moves possibly less than a full Markov basis. If $I$ is an ideal in the ring $Q[x_1,\ldots,x_d]$ and $f$ is a polynomial, then the saturation of $I$ by $f$ (denoted $I:f^\infty$) is defined by

$$I:f^\infty := \{g \in Q[x_1,\ldots,x_d] : f^k \cdot g \in I \text{ for some } k \geq 0\},$$

which is also an ideal. For the indeterminate $x_i$, $I:x_i^\infty$ is the collection of polynomials $g$ such that $x_i^k g$ is in the ideal $I$ for some choice of the exponent $k$.

PROPOSITION 4.1. *Suppose $M_{A,\mathbf{t}}$ is a Markov subbasis, let $M_{A,\mathbf{t}} = \{\pm\mathbf{m}_1, \ldots, \pm\mathbf{m}_g\}$ and let $G := \{\mathbf{x}^{\mathbf{m}_i^+} - \mathbf{x}^{\mathbf{m}_i^-}, i = 1,\ldots,g\}$. Suppose $G$ has the following three properties:* (1) *$G$ is a lex Gröbner basis for the generated ideal $I_{M_{A,\mathbf{t}}}$ with order $x_1 > x_2 > \cdots > x_d$ on indeterminates;* (2) *$G \cap Q[x_i,\ldots,x_d]$ are square-free in $x_i$ for each $i$; and* (3) *$(I_{M_{A,\mathbf{t}}}:x_i^\infty) \cap Q[x_{i+1},\ldots,x_d] \subset I_{M_{A,\mathbf{t}}}$ for each $i=1,2,\ldots,d-1$. Then the polytope $A^{-1}[\mathbf{t}]$ has the sequential interval property.*

PROOF. By Lemma 4.1, $\pi_1(A^{-1}[\mathbf{t}])$ is an interval. We must show that two tables in $A^{-1}[\mathbf{t}]$ with a common entry in coordinate 1 can be connected with moves in $M_{A,\mathbf{t}}$ without touching coordinate 1. To see this, suppose tables $\mathbf{u}', \mathbf{v}' \in A^{-1}[\mathbf{t}]$ have common first coordinate $u_1 = v_1 = c$.

Let $\mathbf{u} = (0, u_2, u_3, \ldots, u_d)$, $\mathbf{v} = (0, v_2, v_3, \ldots, v_d)$. We must show that $\mathbf{x}^{\mathbf{u}} - \mathbf{x}^{\mathbf{v}} \in \langle G \cap Q[x_2,\ldots,x_d]\rangle$ to be able to connect them with moves in $G$ that only involve changing the second coordinate (by only $\pm 1$ at each step). Since $G$ is a lex basis, $\langle G \cap Q[x_2,\ldots,x_d]\rangle = I_{M_{A,\mathbf{t}}} \cap Q[x_2,\ldots,x_d]$, and it is enough to show that $\mathbf{x}^{\mathbf{u}} - \mathbf{x}^{\mathbf{v}} \in I_{M_{A,\mathbf{t}}}$. We have that $\mathbf{x}^{\mathbf{u}} - \mathbf{x}^{\mathbf{v}} \in I_A$.

Since $x_1^c(\mathbf{x}^{\mathbf{u}} - \mathbf{x}^{\mathbf{v}}) = \mathbf{x}^{\mathbf{u}'} - \mathbf{x}^{\mathbf{v}'} \in I_{M_{A,\mathbf{t}}}$, the binomial $\mathbf{x}^{\mathbf{u}} - \mathbf{x}^{\mathbf{v}} \in (I_{M_{A,\mathbf{t}}}:x_1^\infty) \cap Q[x_2,\ldots,x_d]$. Under the assumption $(I_{M_{A,\mathbf{t}}}:x_1^\infty) \cap Q[x_2,\ldots,x_d] \subset I_{M_{A,\mathbf{t}}}$, the first step is proven.

Suppose now that two tables $\mathbf{u}', \mathbf{v}' \in A^{-1}[\mathbf{t}]$ have common first two coordinates $u_1 = v_1 = c_1$, $u_2 = v_2 = c_2$. Let $\mathbf{u} = (0, 0, u_3, u_4, \ldots, u_d)$, $\mathbf{v} = (0, 0, v_3, v_4, \ldots, v_d)$. We must show that $\mathbf{x}^{\mathbf{u}} - \mathbf{x}^{\mathbf{v}} \in \langle G \cap Q[x_3, \ldots, x_d] \rangle$ to be able to connect them with moves in $G$ that only involve changing the third coordinate (by only $\pm 1$ at each step). By the argument above, we have that $x_2^{c_2} \mathbf{x}^{\mathbf{u}} - x_2^{c_2} \mathbf{x}^{\mathbf{v}} \in I_{M_{A,\mathbf{t}}}$. Then by the saturation condition on $x_2$, $\mathbf{x}^{\mathbf{u}} - \mathbf{x}^{\mathbf{v}} \in (I_{M_{A,\mathbf{t}}} : x_2^{\infty}) \cap Q[x_3, \ldots, x_d] \subset M_{A,\mathbf{t}}$.

The argument continues likewise for each cell in the order $1, 2, \ldots, d$. $\square$

To use Proposition 4.1, one must have in hand a Markov subbasis, which requires knowing some connectivity properties. These can be established sometimes with ad hoc arguments or with the primary decomposition of the ideal $I_{M_{A,\mathbf{t}}}$. Lemma 4.2 below is a new method to verify a Markov subbasis, and we use it in Example 7.3. The quotient ":" operation is defined by $I : f := \{g : f \cdot g \in I\}$, the result of one step of the saturation procedure defined above.

LEMMA 4.2. *Let $M \subset \ker_Z(A)$ be Markov moves with ideal $I_M$. Suppose each element $\mathbf{n} \in A^{-1}[\mathbf{t}]$ satisfies $n_s > 0$ for all $s \in S \subset \{1, \ldots, d\}$, and suppose that $(I_M : \prod_{s \in S} x_s) = I_A$, the toric ideal. Then the moves in $M$ connect all of $A^{-1}[\mathbf{t}]$ and are therefore a Markov subbasis.*

PROOF. Let $\mathbf{u}, \mathbf{v} \in A^{-1}[\mathbf{t}]$, and let $\mathbf{u}' = \mathbf{u} - I_S$, $\mathbf{v}' = \mathbf{v} - I_S$, where $I_S$ is the vector with 1 in the coordinates that are in the set $S$, and 0 elsewhere. Clearly, $\mathbf{x}^{\mathbf{u}'} - \mathbf{x}^{\mathbf{v}'} \in I_A$, so by the saturation assumption $(\mathbf{x}^{\mathbf{u}'} - \mathbf{x}^{\mathbf{v}'}) \prod_{s \in S} x_s \in I_M$. The fundamental result of Diaconis and Sturmfels ([14], Theorem 3.1) says that the moves in $M$ connect $\mathbf{u} = \mathbf{u}' + I_S$ with $\mathbf{v} = \mathbf{v}' + I_S$ through the nonnegative tables. $\square$

**5. Bounds on cell entries.** When the conditions for sequential interval property are met, the next question is how to quickly determine or approximate the upper and lower bounds of the interval $[l, u]$. In very special cases one can use known formulas for the interval, such as the Fréchet bounds. This works in two-way tables and some decomposable graphical models [6]. For general multiway tables, usually no simple formula is available to compute the bounds. Three general ways to determine or approximate the upper and lower bounds of the interval $[l, u]$ are integer programming (IP), linear programming (LP) and the shuttle algorithm. IP always gives the exact integer bounds $l$ and $u$, but it is much slower than the other two methods.

LP in the rational numbers can dynamically find bounds on the interval at each step in the sampling. LP is much faster than IP, and under conditions that hold in many examples, LP gives the same answer as IP. The conditions we formulate are concrete algebraic conditions that can be checked with a preliminary calculation. Hosten and Sturmfels [22] study the difference between LP and IP from a different point of view. They give the largest possible difference over all constraint values, whereas our results use the particular constraint values of the data set.

The numerical implementation of LP to determine an interval $[l, u]$ must be done carefully. LP sometimes gives wider intervals than the true interval because LP considers solutions in a larger space. Roundoff of numerical approximations that come from floating point operations or interior point methods can result in sampling a number out of the feasible range $[l, u]$ or into a strict subset of the feasible range which can lead to errors. The program that we embedded into the sampling code and that worked well is lpSolve [1].

A third way to approximate the intervals is the shuttle algorithm, described in [5] and [16]. This is an iterative method that usually does not give exact IP results, but it has two advantages in special cases: it is fast and easy to program, and it can be implemented without explicitly constructing a constraint matrix, a task which may be impossible for very large problems with millions of cells. In our numerical examples LP works better than the shuttle algorithm, in some cases much better.

Consider the IP and LP problems

$$u_j(\mathbf{b}) := \max\{n_j : A_j \mathbf{n} = \mathbf{b}, \mathbf{n} \in Z_+^d\},$$
$$l_j(\mathbf{b}) := \min\{n_j : A_j \mathbf{n} = \mathbf{b}, \mathbf{n} \in Z_+^d\},$$
$$U_j(\mathbf{b}) := \max\{q_j : A_j \mathbf{q} = \mathbf{b}, \mathbf{q} \in Q_+^d\},$$
$$L_j(\mathbf{b}) := \min\{q_j : A_j \mathbf{q} = \mathbf{b}, \mathbf{q} \in Q_+^d\},$$

where $Z_+, Q_+$ are the nonnegative integers and nonnegative rational numbers, respectively. We are interested in bounding the nonnegative quantities $U_j - u_j$ and $l_j - L_j$.

In Propositions 5.1 and 5.2 that follow, we use the relationship between lower and upper IP bounds and normal forms with respect to lex and grevlex term orders explained in [9] and stated in Algorithm 5.6 of [28], page 43. For the following proposition, let $A_Q^{-1}[\mathbf{t}] := \{\mathbf{q} \in Q_+^d : A\mathbf{q} = \mathbf{t}\}$, the set of nonnegative rational vectors with constraints $\mathbf{t}$.

PROPOSITION 5.1. *Suppose a Markov subbasis $M_{A,\mathbf{t}} = \{\pm\mathbf{m}_1, \ldots, \pm\mathbf{m}_g\}$ has the property that $G := \{\mathbf{x}^{\mathbf{m}_i^+} - \mathbf{x}^{\mathbf{m}_i^-}, i = 1, \ldots, g\}$ is a lex Gröbner basis with ordering $x_1 > x_2 > \cdots > x_d$ on indeterminates for the generated ideal $I_{M_{A,\mathbf{t}}}$. Also, suppose $I_{M_{A,\mathbf{t}}} : \prod_{s \in I_{S_Q}} x_s = I_A$, where $S_Q$ is the collection of coordinates which are always positive for elements in $A_Q^{-1}[\mathbf{t}]$, and suppose $I_{M_{A,\mathbf{t}}} : x_i^\infty \cap Q[x_{i+1}, \ldots, x_d] \subset I_{M_{A,\mathbf{t}}}$ for each $i = 1, 2, \ldots, d-1$.*

*If the coordinate values of all* $\mathbf{m}_i^+$ *$(i = 1, \ldots, g)$ are in $\{0, 1\}$, then $l_j(\mathbf{t}_j) = L_j(\mathbf{t}_j)$ for all $j = 1, 2, \ldots, d$ and all $\mathbf{t}_j$ given by $\mathbf{t}_1 = \mathbf{t}$, $\mathbf{t}_j = \mathbf{t} - \mathbf{a}_1 n_1 - \mathbf{a}_2 n_2 - \cdots - \mathbf{a}_{j-1} n_{j-1}$, $j = 2, \ldots, d$.*

PROOF. We show first the result that $l_1 \leq L_1$. Let $\mathbf{m} \in A^{-1}[\mathbf{t}]$.

Use long division to compute the normal form of $\mathbf{x^m}$ with respect to $I_{M_{A,\mathbf{t}}}$. Let the normal form be the monomial $\mathbf{x^{n^\star}}$. It is nearly immediate that $n_1^\star \geq l_1$, since the first coordinate of the normal form when dividing by a Gröbner basis for the full ideal $I_A$ is $l_1$.

Let $\mathbf{q}^\star$ solve $L_1 = \min\{q_1 : A\mathbf{q} = \mathbf{t}, \mathbf{q} \in Q_+^d\}$. We show that $q_1^\star \geq n_1^\star$, which together with $n_1^\star \geq l_1$ will prove the result $L_1 = l_1$.

Suppose by way of contradiction that $n_1^\star > q_1^\star$. Since $\mathbf{q}^\star$ is rational, an integer multiple, say $\lambda \mathbf{q}^\star$, is integral. Then $A(\lambda \mathbf{q}^\star) = A(\lambda \mathbf{n}^\star)$, so $\mathbf{x}^{\lambda \mathbf{n}^\star} - \mathbf{x}^{\lambda \mathbf{q}^\star} \in I_A$. Furthermore, by the assumption of positivity of coordinates $S_Q$ on elements in $A_Q^{-1}[\mathbf{t}]$, it follows that $\mathbf{q}_s^\star, \mathbf{n}_s^\star > 0$ for $s \in S_Q$. Then $\mathbf{x}^{\lambda \mathbf{n}^\star} - \mathbf{x}^{\lambda \mathbf{q}^\star} \in I_{M_{A,\mathbf{t}}}$ by the assumption $I_{M_{A,\mathbf{t}}} : \prod_{s \in I_{S_Q}} x_s = I_A$.

Since $G$ is a Gröbner basis for this ideal, one of the lead terms of the basis must divide the lead monomial $\mathbf{x}^{\lambda \mathbf{n}^\star}$. This means that the indices of positive coordinates of the exponents $\mathbf{m}_i^+$ of the lead monomial must be included in the positive coordinates of $\mathbf{n}^\star$. Since the corresponding coordinate values are 0 or 1, the divisor must also divide $\mathbf{n}^\star$. This contradicts its construction above as the normal form without divisors. Hence, it cannot be the case that $n_1^\star > q_1^\star$. This proves that $l_1 \leq n_1^\star \leq q_1^\star = L_1$.

We show next the result that $l_2 \leq L_2$. Let $\mathbf{m} \in A^{-1}[\mathbf{t}]$.

Use long division to compute the normal form of $\mathbf{x^m}$ with respect to $G_2 := G \cap Q[x_2, x_3, \ldots, x_d]$, the elements of the subbasis that only involve coordinates $2, 3, \ldots, d$. Let the normal form be the monomial $\mathbf{x}^{\mathbf{n}^\star}$, where $n_1^\star = m_1$, which has not changed in the division. It is nearly immediate that $n_2^\star \geq l_2$, since the first coordinate of the normal form when dividing $x_2^{m_2} \cdots x_d^{m_d}$ by a Gröbner basis for the full ideal $I_{A_2}$ is $l_2$.

Let $\mathbf{q}^\star$ solve $L_2 = \min\{q_2 : A\mathbf{q} = \mathbf{t}, q_1 = m_1, \mathbf{q} \in Q_+^d\}$. We show that $q_2^\star \geq n_2^\star$, which together with $n_2^\star \geq l_2$ will prove the result $L_2 = l_2$.

Suppose by way of contradiction that $n_2^\star > q_2^\star$. Since $\mathbf{q}^\star$ is rational, an integer multiple, say $\lambda \mathbf{q}^\star$, is integral. Then $A(\lambda \mathbf{q}^\star) = A(\lambda \mathbf{n}^\star)$, so $\mathbf{x}^{\lambda \mathbf{n}^\star} - \mathbf{x}^{\lambda \mathbf{q}^\star} \in I_A$. Furthermore, by the assumption of positivity of coordinates $S_Q$ on elements in $A_Q^{-1}[\mathbf{t}]$, it follows that $\mathbf{q}_s^\star, \mathbf{n}_s^\star > 0$ for $s \in S_Q$. Then $\mathbf{x}^{\lambda \mathbf{n}^\star} - \mathbf{x}^{\lambda \mathbf{q}^\star} \in I_{M_{A,\mathbf{t}}}$, since $I_{M_{A,\mathbf{t}}} : \prod_{s \in I_{S_Q}} x_s = I_A$. Also, $\mathbf{x}^{\lambda(0, n_2^\star, \ldots, n_d^\star)} - \mathbf{x}^{\lambda(0, q_2^\star, \ldots, q_d^\star)} \in I_{M_{A,\mathbf{t}}} : x_1^\infty \cap Q[x_2, \ldots, x_d]$. By the assumption on $I_{M_{A,\mathbf{t}}} : x_i^\infty$ it follows that $\mathbf{x}^{\lambda(0, n_2^\star, \ldots, n_d^\star)} - \mathbf{x}^{\lambda(0, q_2^\star, \ldots, q_d^\star)} \in I_{M_{A,\mathbf{t}}}$.

Since $G_2$ is a lex Gröbner basis for the ideal $I_{M_{A,\mathbf{t}}} \cap Q[x_2, \ldots, x_d]$ by the elimination theorem, one of the lead terms of the basis $G_2$ must divide the lead monomial $\mathbf{x}^{\lambda(0,n_2^\star,\ldots,n_d^\star)}$, since we have just shown that this is the lead monomial in a binomial that belongs to $I_{M_{A,\mathbf{t}}} \cap Q[x_2,\ldots,x_d]$. This means that the indices of positive coordinates of the exponents $\mathbf{m}_i^+$ of the lead monomial must be included in the positive coordinates of $\mathbf{n}^\star$. Since the corresponding coordinate values are 0 or 1, the divisor must also divide $\mathbf{n}^\star$. This contradicts its construction above as the normal form without divisors. Hence, it cannot be the case that $n_2^\star > q_2^\star$. This proves that $l_2 \leq n_2^\star \leq q_2^\star = L_2$.

The remaining coordinates are proved similarly. □

There is a corresponding result for the upper bounds. Whereas the lex basis relates IP minimization to the normal form of a monomial, it is the grevlex basis that relates IP maximization to the normal form. We state the result below only for the first entry, since it must be applied repeatedly. Using the result requires recomputing a grevlex basis for each of the matrices $A_i$ (containing columns $i, i+1, \ldots, d$ from $A$) and rechecking the condition, because we cannot simply apply an elimination theorem on a single lex basis as before.

PROPOSITION 5.2. *Suppose a Markov subbasis $M_{A,\mathbf{t}} = \{\pm \mathbf{m}_1, \ldots, \pm \mathbf{m}_g\}$ has the property that $G := \{\mathbf{x}^{\mathbf{m}_i^+} - \mathbf{x}^{\mathbf{m}_i^-}, i = 1, \ldots, g\}$ is a grevlex Gröbner basis with ordering $x_d > x_{d-1} > \cdots > x_1$ on indeterminates for the generated ideal $I_{M_{A,\mathbf{t}}}$. Also, suppose $I_{M_{A,\mathbf{t}}} : \prod_{s \in I_{S_Q}} x_s = I_A$, where $S_Q$ is the collection of coordinates which are always positive for elements in $A_Q^{-1}[\mathbf{t}]$. If the coordinate values of $\mathbf{m}_i^+$ are in $\{0,1\}$, then $u_1(\mathbf{t}) = U_1(\mathbf{t})$.*

PROOF. We show that $U_1 \leq u_1$. Let $\mathbf{m} \in A^{-1}[\mathbf{t}]$.

Use long division to compute the normal form of $\mathbf{x}^\mathbf{m}$ with respect to the grevlex basis $I_{M_{A,\mathbf{t}}}$. Let the normal form be the monomial $\mathbf{x}^{\mathbf{n}^\star}$. It is nearly immediate that $n_1^\star \leq u_1$, since the exponent on $x_1$ of the normal form when dividing by a grevlex Gröbner basis with reversed indeterminate order for the full ideal $I_A$ is $u_1$.

Let $\mathbf{q}^\star$ solve $U_1 = \max\{q_1 : A\mathbf{q} = \mathbf{t}, \mathbf{q} \in Q_+^d\}$. We show that $q_1^\star \leq n_1^\star$, which together with $n_1^\star \leq u_1$ will prove the result $U_1 \leq u_1$.

Suppose by way of contradiction that $q_1^\star > n_1^\star$. Since $\mathbf{q}^\star$ is rational, an integer multiple, say $\lambda \mathbf{q}^\star$, is integral. Then $A(\lambda \mathbf{q}^\star) = A(\lambda \mathbf{n}^\star)$, so $\mathbf{x}^{\lambda \mathbf{n}^\star} - \mathbf{x}^{\lambda \mathbf{q}^\star} \in I_A$. Furthermore, by the assumption of positivity of coordinates $S_Q$ on elements in $A_Q^{-1}[\mathbf{t}]$, it follows that $\mathbf{q}_s^\star, \mathbf{n}_s^\star > 0$ for $s \in S_Q$. Then $\mathbf{x}^{\lambda \mathbf{n}^\star} - \mathbf{x}^{\lambda \mathbf{q}^\star} \in I_{M_{A,\mathbf{t}}}$ by the assumption $I_{M_{A,\mathbf{t}}} : \prod_{s \in I_{S_Q}} x_s = I_A$.

Since $G$ is a Gröbner basis for this ideal, one of the lead terms of the basis must divide the lead monomial $\mathbf{x}^{\lambda \mathbf{n}^\star}$. This means that the indices of positive coordinates of the exponents $\mathbf{m}_i^+$ of the lead monomial must be included in the positive coordinates of $\mathbf{n}^\star$. Since the corresponding coordinate values are 0 or 1, the divisor must also divide $\mathbf{n}^\star$. This contradicts its construction above as the normal form without divisors. Hence it cannot be the case that $q_1^\star > n_1^\star$. This proves that $U_1 = q_1^\star \leq n_1^\star \leq u_1$. $\square$

The corollary below, combining Propositions 5.1 and 5.2, applies directly to Examples 7.1, 7.2 and 7.4.

COROLLARY 5.1. *If a lex Gröbner basis for $I_A$ has square-free exponents on the lead monomials, then $l_j = L_j$ for all $j = 1, \ldots, d$. If each grevlex Gröbner basis for $I_{A_j}, j = 1, \ldots, d$, and indeterminate ordering $x_d > x_{d-1} > \cdots > x_1$ has square-free exponents on the lead monomials, then $u_j = U_j$ for all $j = 1, \ldots, d$.*

PROOF. The assumptions of Proposition 5.1 hold if $I_{M_{A,\mathbf{t}}} = I_A$, so the lower bounds from LP and IP are equal. For the upper bounds, the statement is a restatement of Proposition 5.2 for each step in the sequential sampling. $\square$

**6. Sampling distributions.** Assume that the sequential interval property holds for a multiway table with given constraints, and that the intervals can be approximated by LP. The next question is how to sample from these intervals. Ideally, we want to sample a cell value from the true marginal distribution of a cell entry conditional on the entries that have already been sampled. However, these marginal distributions are quite difficult to compute explicitly except in very small examples. SIS samples from a simple proposal distribution (rather than the true distribution) on the set of all possible marginal values.

For a target uniform distribution, which is useful for counting the total number of tables and some Bayesian applications, we propose a uniform distribution on the available interval for each cell, that is, $p(x) = 1/(u-l+1)$ on integers in the interval $[l, u]$. We call this the "uniform sampling method." With the length of the proposed sampling interval, the importance weights can be computed exactly for reweighting at the end. This strategy gives low $cv^2$ ($\leq 5$ for all examples we have tested) and works very well on the examples in Section 7.

For a target hypergeometric distribution, which arises in conditional inference with multinomial sampling, we propose to sample a cell value from the hypergeometric distribution $p(x) = \binom{u}{x}\binom{u}{l+u-x}/\binom{2u}{l+u}$ on the interval of available integers $[l, u]$. We call this the "hypergeometric sampling method," which is usually (but not always, see Example 7.4) better than the uniform sampling method when the target distribution is hypergeometric. This hypergeometric proposal does not give the exact hypergeometric target in the end. It is just a reasonable marginal approximation. This method gives satisfactory results for examples in Section 7, although the $cv^2$ is not consistently small. For sparse tables, approximating the marginal mass function of the count in a single cell can be difficult.

**7. Examples.** In the examples that follow, we sample sequentially from intervals computed with the LP approximation. The LP approximation is very close to or exactly equal to the IP range in all examples. In Example 7.1 one can apply known results on Markov bases to avoid algebraic calculations, and the most basic results of Section 3 apply. In Example 7.2 one must do explicit algebraic calculations to verify the conditions of Section 3. We did a detailed numerical comparison with the Markov chain on Example 7.2. Example 7.3 (6-way Czech autoworker data) is one that requires the full theory of Markov subbases of Section 4 and consideration of the specific margin values to get sequential intervals under one model. We also study a second model for which we could not compute the Markov basis, and we see that SIS still works well. The no-3-way interaction model of Example 7.4 is a well-known example where the Gröbner basis involves moves of size 2, and yet the sequential theory applies perfectly. Example 7.5 is a classic triangular genotype table, and it brings out the importance of checking different cell orders. In some orders the sequential interval property holds, and in other quite natural orders it does not, and this can be seen in the lex basis. Finally, Example 7.6 is an important application of sampling on lattice points that are not strictly speaking contingency tables. The work of Rapallo [27] on Markov bases and structural zeros may be useful for other examples.

The starting point to verify the conditions of Sections 3, 4 and 5 for a particular example is to attempt to compute the toric ideal $I_A$. For this we have used the toric library `toric.lib` in the free software Singular [19] and the `groebner` command in $4ti2$ [21]. The software $4ti2$ was used to construct constraint matrices for several examples. The operations of saturation and quotient (":") that figure in the results of Sections 4 and 5 were done quickly in Singular.

In the following examples, all results are based on 1000 random samples using either the uniform sampling method or the hypergeometric sampling method. The code was written in R [26] and the software lpSolve was called from R. The running times range from several seconds to a few minutes on a 2.0 GHz computer. When IP is used instead of LP, a computation typically takes hours, and sometimes it will not terminate in a reasonable amount of time.

EXAMPLE 7.1. Consider the 3-way case/control data (Table 1) in the $4 \times 4 \times 2$ table from the Ille-et-Verlaine cancer study of the age 35–44 group

SIS FOR MULTIWAY TABLES 17

Table 1
*Age* 35–44 *data on oesophageal cancer from* [3]

|       |   |   | \ **A** | | | |
|-------|---|---|----|----|----|----|
|       |   |   | 1 | 2 | 3 | 4 |
| R = 0 | T | 1 | 60 | 35 | 11 | 1 |
|       |   | 2 | 13 | 20 | 6  | 3 |
|       |   | 3 | 7  | 13 | 2  | 2 |
|       |   | 4 | 8  | 8  | 1  | 0 |
| R = 1 | T | 1 | 0  | 0  | 0  | 2 |
|       |   | 2 | 1  | 3  | 0  | 0 |
|       |   | 3 | 0  | 1  | 0  | 2 |
|       |   | 4 | 0  | 0  | 0  | 0 |

([3], Appendix I). The factors are Alcohol level (A), Tobacco level (T) and Response R, where R = 0 is a control measurement and R = 1 is a case.

The "case" outcomes are sampled with a multinomial distribution with probabilities $p(a,t|1)$ on the Alcohol and Tobacco covariates. The "control" outcomes are also sampled with a multinomial distribution with probabilities $p(a,t|0)$. With a retrospective model $p(a,t|1)/p(a,t|0) = e^{\alpha_a + \beta_t}$ of eight parameters, the appropriate margins to fix for conditional inference [treating $p(a,t|0)$ as unknown nuisance parameters] are [A, T] (sum over case/control counts at each level), [A, R] and [T, R] (sums over other factor at each response level). The constraints imply that the Graver basis ([28], page 55) for the independence model on T and R is a Markov basis, and the Graver basis is equivalent to the collection of square-free circuit moves on one level of the Response factor. Thus, the results of Section 3 and Corollary 5.1 imply the property of sequential intervals and LP will give the exact integral interval bounds at each step.

The simulation with LP gave 100% good tables. When the underlying distribution is uniform, the uniform sampling method gave $cv^2$ of 0.24 and estimated the total number of tables to be 25, a number confirmed by LattE in a total elapsed time of 7 seconds on a 2.8 GHz desktop. When the underlying distribution is hypergeometric, the hypergeometric sampling method gave $cv^2$ of 0.5, and the estimated $p$-value for the exact goodness-of-fit test [defined by equations (1) and (2)] is 0.04.

EXAMPLE 7.2. Consider the 4-way abortion opinion data (Table 2) from [8], page 129. The observations are classified according to race, sex, age and opinion. There are three different opinions: yes means supporting legalized abortion, no means opposing legalized abortion, and the last one is undecided.

TABLE 2
*4-way abortion opinion data from* [8]

| Race | Sex | Opinion | 18–25 | 26–35 | 36–45 | 46–55 | 56–65 | 66+ |
|---|---|---|---|---|---|---|---|---|
| White | Male | Yes | 96 | 138 | 117 | 75 | 72 | 83 |
| | | No | 44 | 64 | 56 | 48 | 49 | 60 |
| | | Undec. | 1 | 2 | 6 | 5 | 6 | 8 |
| | Female | Yes | 140 | 171 | 152 | 101 | 102 | 111 |
| | | No | 43 | 65 | 58 | 51 | 58 | 67 |
| | | Undec. | 1 | 4 | 9 | 9 | 10 | 16 |
| Nonwhite | Male | Yes | 24 | 18 | 16 | 12 | 6 | 4 |
| | | No | 5 | 7 | 7 | 6 | 8 | 10 |
| | | Undec. | 2 | 1 | 3 | 4 | 3 | 4 |
| | Female | Yes | 21 | 25 | 20 | 17 | 14 | 13 |
| | | No | 4 | 6 | 5 | 5 | 5 | 5 |
| | | Undec. | 1 | 2 | 1 | 1 | 1 | 1 |

Christensen fits the log-linear model for the expected cell counts with all three-way interactions and all lower order terms. A shorthand notation for this model is to list its highest-order interaction terms: [RSO], [RSA], [ROA] and [SOA]. The conditional goodness-of-fit test for this model requires fixing all 3-way margins, [R, S, O], [R, S, A], [R, O, A] and [S, O, A]. The lex basis of 165 elements is square-free in the lead monomials, so the sequential interval property holds by Section 3 and the IP and LP lower bounds are identical. A more detailed calculation to verify the conditions of Corollary 5.1 requires computing a grevlex basis for each of the submatrices of $A_i$, defined in Section 3 as the matrix that has columns $i, i+1, \ldots, d$ from $A$. This can be done and the condition is verified, proving that LP and IP upper bounds are always the same.

The LP method for finding the interval bounds gave 100% good tables in practice. When the underlying distribution is uniform, the uniform sampling method gave $cv^2$ of 2.92 and estimated the total number of tables to be $9.1 \times 10^7$. When the underlying distribution is hypergeometric, the value of $cv^2$ using the hypergeometric sampling method was around 102.9, and the estimated $p$-value for the exact goodness-of-fit test [defined by (1) and (2)] is 0.85 with standard error 0.1, based on 1000 tables which took about 5 minutes in R on a 2.0 GHz computer. The MCMC algorithm generated 1000 samples (with 1,000,000 samples as burn-in) in 224 minutes and estimated the $p$-value to be 0.84 with standard error 0.05. Thus, SIS is about 11 times faster than the MCMC algorithm for this example.

The algebraic conditions for SIS with some models on this data are difficult to verify. For example, 4*ti*2 runs for an hour on a 2.8 GHz Linux desktop with 1 GB of memory without completing the Markov basis calculation on the model [RS], [RA], [RO], [SO], [SA], [OA].

TABLE 3
*6-way Czech autoworker data from* [18]

| | | | | | B | no | | yes | |
|---|---|---|---|---|---|---|---|---|---|
| F | E | D | C | A | no | yes | no | yes |
| Negative | <3 | <140 | no | | 44 | 40 | 112 | 67 |
| | | | yes | | 129 | 145 | 12 | 23 |
| | | ≥140 | no | | 35 | 12 | 80 | 33 |
| | | | yes | | 109 | 67 | 7 | 9 |
| | ≥3 | <140 | no | | 23 | 32 | 70 | 66 |
| | | | yes | | 50 | 80 | (0) 7 | 13 |
| | | ≥140 | no | | 24 | 25 | 73 | 57 |
| | | | yes | | 51 | 63 | 7 | 16 |
| Positive | <3 | <140 | no | | 5 | 7 | 21 | 9 |
| | | | yes | | (0) 9 | 17 | (0) 1 | (0) 4 |
| | | ≥140 | no | | (0) 4 | 3 | 11 | 8 |
| | | | yes | | 14 | 17 | 5 | (0) 2 |
| | ≥3 | <140 | no | | 7 | (0) 3 | 14 | 14 |
| | | | yes | | 9 | 16 | (0) 2 | (0) 3 |
| | | ≥140 | no | | (0) 4 | (0) 0 | 13 | 11 |
| | | | yes | | (0) 5 | 14 | (0) 4 | 4 |

EXAMPLE 7.3. Consider the 6-way binary Czech autoworker data in Table 3 from a prospective study of probable risk factors for coronary thrombosis [18]. There are 1,841 men in a car factory involved in the study. Here A, B, C, D, E and F indicate different risk factors. One reasonable model is given by [ACDEF], [ABDEF], [ABCDE], [BCDF], [ABCF], [BCEF] [17]. The conditional goodness-of-fit test for this model requires fixing the three 5-way and the three 4-way margins in the above model representation. Implementing SIS for this example requires techniques beyond the basic methods of Section 3, because the lex basis does not have square-free lead exponents.

In Table 3 (0) indicates that the LP lower bound for that cell entry is 0 with the constraints from the model above; the others are strictly positive. Identifying these cells is relevant when we apply Propositions 4.1, 5.1 and 5.2, as the (0) cells form the complement of the set $S_Q$ (defined in Proposition 5.1).

The lex basis for the toric ideal with lex order in indeterminates yields 20 elements, the first of which has an exponent of 2 on the lead indeterminate $x_{111111}$. Therefore, Proposition 3.1 cannot be applied directly. However, the ideal generated by the other 19 polynomials saturates in one step with respect to the monomial $\prod_{s \in S} x_s$, where $S$ is the set of 41 coordinates that must be positive. Hence, by Lemma 4.2 these 19 moves are a Markov subbasis. They are a lex Gröbner basis for themselves, and they have the saturation property required in Proposition 4.1, so the sequential interval property holds. Furthermore, Proposition 5.1 shows that the IP and LP lower bounds are always the same [which also implies that the (0) cells in the rationals are the same cells as those that could be 0 in the integers]. Corollary 5.1 does not apply to show that the LP and IP upper bounds are the same, because exponents of 2 appear in the grevlex bases. We can use Proposition 5.2 on successive cells to show that LP and IP are the same after a few initial cells.

If the cells are filled in across rows and then down, the order is $111111, 211111, 121111,\ldots$ (the order from $4ti2$). The sequential interval property holds in this order as well.

For this model, using LP for interval bounds gave 100% good tables. The shuttle algorithm gave 99% good tables with one iteration and 99.5% with two iterations. When the underlying distribution is uniform, the uniform sampling method gave $cv^2$ of 1.09 and estimated the total number of tables to be 841. The quantity $cv^2$ when targeting the hypergeometric distribution using the hypergeometric sampling method was 50.7, and the estimated $p$-value for the exact goodness-of-fit test [defined by (1) and (2)] is 0.27. Fitting this model in R using the `loglin` command gives a $\chi^2$ statistic of 5.8 on 4 degrees of freedom, for a $p$-value of approximately 0.21.

Consider the model of all 15 four-element constraints like $[A, B, C, D]$, that is, all 4-way margins. We could not obtain the Markov basis for this model, but SIS still works well with $cv^2 = 5.0$ when the target distribution is uniform. LP gave 100% good tables, whereas the shuttle algorithm gave only 2% good tables after 10 iterations.

EXAMPLE 7.4. Consider the $3 \times 3 \times 3$ example (Table 4) from [14], page 379, with a model of no-3-way interaction. The conditional goodness-of-fit test for this model requires fixing all "line sums."

When ordered left to right across rows, Proposition 3.1 implies sequential intervals and Corollary 5.1 gives an IP/LP gap of 0 at every step. In simulation LP gave 100% good tables, and the shuttle algorithm also gave 100% good tables after one iteration.

When the underlying distribution is uniform, the uniform sampling method gave $cv^2$ of 2.08 and estimated the total number of tables to be $1.9 \times 10^{12}$. This is consistent with the number 1,919,899,782,953 from LattE, computed in a total elapsed real time of 45 seconds on a 2.8 GHz desktop computer.

TABLE 4
$3 \times 3 \times 3$ table from [14]

| 9  | 16 | 41  | 8  | 8  | 46 | 11 | 14 | 38  |
|----|----|-----|----|----|----|----|----|-----|
| 85 | 52 | 105 | 35 | 29 | 54 | 47 | 35 | 115 |
| 77 | 30 | 38  | 37 | 15 | 22 | 25 | 21 | 42  |

Table 5
*Rhesus data from [20]*

|       | $B_1$ | $B_2$ | $B_3$ | $B_4$ | $B_5$ | $B_6$ | $B_7$ | $B_8$ | $B_9$ |
|-------|------|-----|----|------|----|------|----|----|----|
| $B_1$ | 1236 |     |    |      |    |      |    |    |    |
| $B_2$ | 120  | 3   |    |      |    |      |    |    |    |
| $B_3$ | 18   | 0   | 0  |      |    |      |    |    |    |
| $B_4$ | 982  | 55  | 7  | 249  |    |      |    |    |    |
| $B_5$ | 32   | 1   | 0  | 12   | 0  |      |    |    |    |
| $B_6$ | 2582 | 132 | 20 | 1162 | 29 | 1312 |    |    |    |
| $B_7$ | 6    | 0   | 0  | 4    | 0  | 4    | 0  |    |    |
| $B_8$ | 2    | 0   | 0  | 0    | 0  | 0    | 0  | 0  |    |
| $B_9$ | 115  | 5   | 2  | 53   | 1  | 149  | 0  | 0  | 4  |

When targeting the hypergeometric distribution, the hypergeometric sampling method gave $cv^2 = 180.7$.

EXAMPLE 7.5. Consider data in Table 5 of genotype pairs from [20]. The constraints for conditional goodness-of-fit test of Hardy–Weinberg proportions are the nine allele counts, which are nine linear functions that count twice the diagonal entry, so the $A$ matrix has entries 0, 1 and 2. For sequential sampling, the order of cells given by Table 6 leads to sequential intervals by Proposition 3.1. In general, for the genotype problem sampling across rows will not give intervals.

The lead monomials in a lex basis have exponents that are all 0 or 1, so LP gives the exact lower bounds by Proposition 5.1. For the upper bound, the grevlex condition of Proposition 5.2 does not hold from the first cell, but it does hold after a few cells, so IP and LP give the same bounds after some initial cells. The simulation with LP produced 100% good tables. See [23] for a direct sampling strategy and some further discussion of this example.

EXAMPLE 7.6. Consider a constraint matrix $A$ of the form $A = (A_0|I)$ with 0 or 1 entries. Here $I$ is the $e \times e$ identity matrix and $A_0$ is size $e \times f$ with columns $\mathbf{a}_1,\ldots,\mathbf{a}_f$. This occurs in a tomography problem introduced by Vardi [31], where $A$ is a routing matrix for which routes between adjacent vertices use the connecting edge, and the edge counts are put last as slack variables. The integer data $\mathbf{y} = A\mathbf{x}$, where $\mathbf{x}$ are traffic counts between ordered pairs of nodes on a graph and $\mathbf{y}$ is the aggregate traffic across links. The sampling method of Tebaldi and West [30] for Bayesian computation of the posterior distribution is closely related to sequential sampling. Dinwoodie [15] shows how fast sampling can be used in a Monte Carlo EM algorithm for estimating traffic rates.

Table 6
*Order of cells*

| 1  |    |    |    |    |    |    |    |   |
|----|----|----|----|----|----|----|----|---|
| 10 | 2  |    |    |    |    |    |    |   |
| 11 | 18 | 3  |    |    |    |    |    |   |
| 12 | 19 | 25 | 4  |    |    |    |    |   |
| 13 | 20 | 26 | 31 | 5  |    |    |    |   |
| 14 | 21 | 27 | 32 | 36 | 6  |    |    |   |
| 15 | 22 | 28 | 33 | 37 | 40 | 7  |    |   |
| 16 | 23 | 29 | 34 | 38 | 41 | 43 | 8  |   |
| 17 | 24 | 30 | 35 | 39 | 42 | 44 | 45 | 9 |

The property of sequential intervals holds for the entries of $\mathbf{x}$ under the constraint $A\mathbf{x} = \mathbf{y}$ in the order of the columns. With indeterminates $w_1,\ldots,w_f$ for the first $f$ columns and $z_1,\ldots,z_e$ for the last $e$ slack variables, a lex Gröbner basis in $Q[w_1,\ldots,w_f,z_1,\ldots,z_e]$ consists of the $f$ binomials $w_i - \mathbf{z}^{\mathbf{a}_i}$.

Linear programming will give the exact interval bounds at each step because of the square-free lead monomials. The shuttle algorithm will also give the exact intervals in one step. The interval for the first cell is exactly $[0, \min_{\{i\colon a_{i,1}>0, 1\leq i\leq f\}}\{y_i\}]$ and the same type of problem recurs at each step $1,\ldots,f$.

It is possible to establish properties of SIS for some classes of examples, or, in other words, for some types of contingency tables with certain constraints and margin values. This is an area of ongoing work, but at this time we can make some statements. Logistic regression tables with one integer covariate and positive column sums (at least one measurement at each level of the covariate) have the sequential interval property. This is proved in [7]. The subbasis that corresponds to differences of adjacent minors satisfies the conditions of Proposition 4.1. However, the IP/LP gap may not be zero.

Also, two-way tables with structural zeros and fixed row and column sums have sequential intervals. The same algebraic technology also shows that case/control data with two factors, such as Example 7.1, has the sequential interval property. We conjecture that decomposable graphical models will have the sequential interval property under some order on the cells, but at this time a careful proof is not complete.

**8. Conclusion.** We have described an efficient sequential importance sampling method for sampling multiway tables with given constraints. It can be used to approximate exact conditional inference on contingency tables. SIS sequentially builds up the proposal distribution by sampling table entries one by one. We have presented a theory that relates algebraic properties of collections of Markov moves to certain geometric properties of contingency tables. The geometric properties of "sequential intervals" and the relationship of IP to LP are important for the performance of sequential sampling.

Many real examples show that the theory is applicable and useful, and can be used in some examples when a Markov basis cannot be found.

In practice, one may try sequential sampling even if the sequential interval property does not hold or if the algebraic conditions are not satisfied or not checked. If one can find rough bounds for each entry and design the proposal distribution carefully, so that the fraction of valid tables is high and the $cv^2$ is low, SIS may still give satisfactory results.

Further work is required to formulate a method to design the proposal distribution at each step. We have seen that the uniform sampling method works very well when the underlying distribution is uniform. However, when the target distribution is hypergeometric, the hypergeometric sampling method could be improved.

**Acknowledgments.** We have used the software 4*ti*2, lpSolve, R and Singular for computations. We thank Sam Buttrey for the lpSolve package in R and Raymond Hemmecke for a version of 4*ti*2. A referee has suggested that Proposition 4.1 holds with the weaker condition $(I_M : \prod_{s \in S} x_s^{\min\{n_s : \mathbf{n} \in A^{-1}(\mathbf{t})\}}) = I_A$ and we think that further theoretical work on connectivity of Markov chains with the goal of computational efficiency would be valuable.

Y. Chen
Department of Statistics
University of Illinois at Urbana-Champaign
Champaign, Illinois 61820
USA
E-mail: yuguo@uiuc.edu

I. H. Dinwoodie
Institute of Statistics
  and Decision Sciences
Duke University
Durham, North Carolina 27708-0251
USA
E-mail: ihd@stat.duke.edu

S. Sullivant
Department of Mathematics
Harvard University
One Oxford Street
Cambridge, Massachusetts 02138
USA
E-mail: seths@math.harvard.edu